# A METHOD WITH FEEDBACK FOR AGGREGATION OF GROUP INCOMPLETE PAIR-WISE COMPARISONS

**Authors:**
**Vitaliy Tsyganok* (ORCID 0000-0002-0821-4877), Sergii Kadenko*, Oleh Andriichuk*, Pavlo Roik***

*Institute for information recording of the National academy of sciences of Ukraine
2, Shpak str.
03113, Kyiv, Ukraine

**ABSTRACT**
A method for aggregation of expert estimates in small groups is proposed. The method is based on combinatorial approach to decomposition of pair-wise comparison matrices and to processing of expert data. It also uses the basic principles of AHP/ANP approaches, such as building of criteria hierarchy to decompose and describe the problem, and evaluation of objects by means of pair-wise comparisons. It allows to derive priorities based on group incomplete pair-wise comparisons and to organize feedback with experts in order to achieve sufficient agreement of their estimates. Double entropy inter-rater index is suggested for usage as agreement measure. Every expert is given an opportunity to use the scale, in which the degree of detail (number of points/grades) most adequately reflects this expert's competence in the issue under consideration, for every single pair comparison. The method takes all conceptual levels of individual expert competence (subject domain, specific problem, individual pair-wise comparison matrix, separate pair-wise comparison) into consideration. The method is intended to be used in the process of strategic planning in weakly-structured subject domains.

## 1. Introduction

During the last several decades expert estimation proved to be the most effective (and sometimes the only) reliable approach, allowing people to make competent and informed decisions in weakly structured subject domains. Weakly-structured domains are domains, influenced by multiple inter-related factors, or criteria, both tangible and intangible. Not all of these factors can be described by quantitave indicators. The only way to describe decision variants according to intangible factors (criteria) is to ask experts from the target subject domain to evaluate them. Aggregation of data, obtained from multiple experts to facilitate decision-making in a given domain, described by multiple criteria, calls for particular techniques. These techniques were developed and improved for decades, resulting in emergence of several areas of decision science, such as group decision-making, multi-criteria decision aid (MCDA), multi-criteria decision-making (MCDM), and others.

Group decision-making represents one of the most important components of decision science, as, in order to improve the reliability of the information obtained from experts it is preferable to use a group of experts, rather than just an individual expert, no matter how competent he or she is (as stated, among others, by (Saaty & Peniwati, 2007)). Beside that, obtaining of data from a group of experts allows the decision-maker to utilize the redundancy of information, which is a powerful and universal property. Various approaches have been utilized to facilitate group decision-making in different areas. Most recent academic efforts in the area of group decision-making were covered, among others, in the following papers: (He, He, & Huang, 2017), (Khaleie & Fasanghari 2012), (Lu, Zhang, & Ruan, 2008), (Mohammadi & Makui 2016), (Pérez, Wikström, Mezei et al., 2013), (Xu & Chen, 2008), (Zhang, Ma, Liu et al. 2012).

The best way to evaluate decision variants, that cannot be numerically described, is to compare them among themselves. This suggestion, made and empirically proved by Tom Saaty and his followers (Saaty, 2008), provided the basis for a whole family of AHP/ANP expert evaluation methods, that is widely used by decision-makers around the world. Pair-wise comparisons still remain in the focus of academic attention. As of now, pair-wise comparison-based approaches are used in combination with neural networks, fuzzy logic, group theory, and even genetic algorithms. Some recent publications on pair-wise comparison-based approaches include (Dong et al. 2015), (Fedrizzi & Brunelli, 2010), (Krejci, 2015), (Zhang & Chen, 2016), (Cavallo & D'Appuzzo, 2012), (Tsyganok, Kadenko, & Andriichuk, 2015).

Aggregation of pair-wise comparison matrices (PCM) is an important problem, for which multiple solution methods were suggested, including geometric mean, arithmetic mean, logarithmic least squares, and others. Well-known approaches to aggregation of individual expert estimates are described by (Forman & Peniwati, 1998). Comparative analysis of different PCM aggregation methods was attempted by many authors, including, for example, (Choo & Wedley, 2004), (Tsyganok, 2010), (Lundi, Siraj, & Greco, 2016).

We feel that publications, dedicated to PCM aggregation methods themselves, and to their comparative analysis, seem to neglect such important issues as expert estimation in groups and verification of consistency (and compatibility) level of expert estimates prior to their aggregation. (Saaty, 1996) suggested that if C.R. was below threshold value, it was necessary to reconsider the problem and revise the judgments. Some efforts in the field of pair-wise comparisons' consistency evaluation and improvement have been made by (Brunelli & Fedrizzi, 2011), (Iida, 2009), (Siraj, Mikhailov, & Keane, 2012 (a)). It is natural to suggest, that consistency and compatibility of expert data should be improved through changes of the most inconsistent ("outstanding") pair-wise comparisons, as it is done in the famous Delphi method. In order to achieve the level of consistency of individual pair-wise comparisons, provided by experts, which would be sufficient for aggregation of these pair-wise comparisons, it is necessary to organize feedback with experts.

In this paper we are going to provide a step-by-step description of an original method of aggregation of group expert pair-wise comparisons with feedback, which has several advantages over other existing approaches. Firstly, it can be used for aggregation of data from both complete and incomplete pair-wise comparison matrices, provided in different scales. Secondly, it thoroughly exploits the redundancy of expert information. Thirdly, it provides an opportunity to organize step-by-step feedback with experts to improve both consistency and compatibility of expert estimates (if necessary) during group decision-making. These and other advantages of the approach will be described in greater detail in the subsequent section of out paper.

## 2. Relevance

Originally, combinatorial approach was suggested by (Tsyganok, 2000). The benefits of combinatorial aggregation method in terms of effectiveness were shown in (Tsyganok, 2010). As our academic school was not the only one working on the combinatory approach, the method was independently re-discovered by (Siraj, Mikhailov, & Keane, 2012 (b)) and respective references were made in (Siraj, Mikhailov, & Keane, 2012 (c)), which also proves the relevance of research. Another brief description of combinatorial algorithm of aggregation was presented in (Tsyganok, Kadenko & Andriichuk, 2015), where the authors primarily focused on the multi-scale approach to expert estimation.

In (Tsyganok, 2010), and (Lundi, Siraj, & Greco, 2016) the method was presented in comparison to other aggregation methods. In (Siraj, Mikhailov, & Keane, 2012 (b)) and (Lundi, Siraj, & Greco, 2016) the authors did not provide any ideas on how to organize feedback with experts if it is necessary for consistency and compatibility improvement (as it was not their purpose). At the same time, in the original publication by (Tsyganok, 2000) the idea of feedback was already present. Beside that, in the listed publications the method was not described as a group one, i.e. it focused on a pair-wise comparison matrix, provided by just one individual expert. In (Tsyganok, Kadenko & Andriichuk, 2015) the group method was described only briefly as part of a strategic planning technology. The technology includes the following conceptual steps:

1) Strategic goal formulation;
2) Selection of experts;
3) Building of criteria hierarchy;
4) Expert estimation of impacts of criteria in the hierarchy;
5) Selection of the optimal strategy of strategic goal achievement. The strategy is the distribution of limited resources among projects, that influence the achievement of strategic goal. This resource allocation allows to achieve the strategic goal with maximum efficiency (at a given moment in time).

In this step-by-step procedure combinatorial algorithm is just a sub-step of step 4). It is utilized to aggregate pair-wise comparisons, provided by a group of experts, and, if necessary, to improve consistency and compatibility of estimates, given by different experts.

As we see, there is no publication, dedicated primarily to combinatorial aggregation method, focusing on its group nature, and describing it in its thoroughness, that would incorporate all the recent research efforts, provide guidelines on how to organize feedback with experts, and that could be used as a universal reference by other authors.

## 3. Problem statement

As it has been mentioned in the previous section, combinatorial algorithm (also called "enumeration of all spanning trees") is intended for aggregation of data, provided by experts in the form of pair-wise comparison matrices. So, in order to approach the aggregation problem in the context of strategic planning technology, described in (Tsyganok, Kadenko & Andriichuk, 2015), we should assume, that:

- we already have a strategic goal, formulated by the decision-maker
- this goal is already decomposed by experts into a hierarchy of factors or criteria (as it is done in AHP/ANP)
- several experts are providing judgments (in the form of PCM), representing pair-wise comparisons of criteria (projects, alternatives)
- we are at a point, where individual expert judgments should be aggregated, so that strategic planning process can move on to its next phase.

With these assumptions in mind, the formal aggregation problem statement looks as follows.

The formal statement of alternative weight calculation problem in our case is as follows.

What is given: $\{A_l$: $l$=[1..$m$]$\}$ – expert pair-wise comparison matrices (PCM) with dimensionality of $n$x$n$, which have the following properties:
1) matrices are reciprocally symmetrical;
2) matrices are multiplicative;
3) in the general case, matrices are incomplete;
4) every single element of a PCM is obtained in some scale, which is assigned a weight coefficient $\{s_h$: $h$=[0..8]$\}$; 0 is the "weight" of a missing element in the PCM, while the largest weight should be assigned to pair-wise comparisons, provided in the scale with the largest number of grades (the most informative one). We suggest calculating scale weights according to Hartley's formula (Hartley, 1928) as explained below.

$c_l$: $l$= [1..$m$]$\}$ – relative competence of experts in the group $\sum_{l=1}^{m} c_l = 1$.

We should find: The resulting object weight (priority) vector $\{w_k$: $k$=[1..$n$]; $\sum_{k=1}^{n} w_k = 1\}$.

## 4. Stages of problem solution

### 4.1 Bringing the estimates given by different experts in different scales to a unified scale.

This step is described in detail in our publication, dedicated to the multi-scale approach to expert estimation (Tsyganok, Kadenko & Andriichuk, 2015). The unified scale is the most detailed (i.e. informative) scale among the scales, used for estimation. This is the scale with the largest number of grades. Judgments, provided in less detailed scales, are approximated by respective grades from the more detailed scale as shown in (Tsyganok, Kadenko & Andriichuk, 2015) and assigned weights, given in the problem statement in section 3 of this paper. We assume that if an expert provides an estimate in the scale with the larger number of grades, this estimate is more informative, and the expert is more confident in his or her knowledge regarding the respective pair-wise comparison.

### 4.2 Weighting of estimates

An estimate provided in the more informative scale should "weight" more than the estimate, provided by an expert using the less detailed scale. The weight should be proportional to the quantity of information, which can be conveyed by means of the respective scale (according to Hartley's formula: $S_N = I = \log_2 N$, where $N$ is the number of expert estimation scale grades (Hartley, 1928)).

If we are dealing with a specific set of pair-wise comparisons (spanning tree) number $q$, reconstructed based on the PCM, provided by expert number $k$, then we suggest calculating the weight coefficient of this set of pair-wise comparisons as geometric mean (GM) of the quantities of information, expressed by the scales, in which all pair-wise comparisons from this basic set were provided (GM of weights of all spanning tree elements):

$$s^{kq} = (\prod_{u=1}^{n-1} \log_2 N_u^{(kq)})^{\frac{1}{n-1}} \quad (1)$$

In this case the number of elements in a spanning tree is $n$ – 1, where $n$ is the number of alternatives (dimensionality of PCM). $N_u^{(kq)}$ is the number of grades in the scale, in which basic pair-wise comparison number $u$ is provided.

If we are dealing with an individual PCM, we, similarly, suggest calculating its weight in terms of the quantity of information as GM of the weights of PCM elements. As PCM are reciprocally symmetrical (according to the problem statement from section 3 of this paper), we can analyze only $n(n-1)/2$ elements, lying above the principal diagonal of the PCM. So, weight of an individual PCM (in terms of the quantity of information), provided by expert number $k$, can be calculated as follows:

$$s^k = (\prod_{\substack{u,v=1 \\ v>u}}^{n} \log_2 N_{uv}^{(k)})^{\frac{2}{n(n-1)}} \quad (2)$$

$N_{uv}^{(k)}$ is the number of grades in the scale, in which the element of the PCM $a_{uv}^{(k)}$ is provided.
number $u$ is provided.
We suggest using GM, because we are dealing with multiplicative scales. However, without loss of generality we could use arithmetic mean as well.

### 4.3 Checking completeness and agreement (consistency and compatibility) of the estimates

Before aggregation of estimates, completeness and agreement of individual PCMs should be checked.

#### 4.3.1 Completeness

Procedure of completeness check is described in detail in (Tsyganok, Kadenko & Andriichuk, 2015). In combinatorial method priority vectors are calculated based on basic sets of pair-wise comparisons (also called spanning trees). If such basic set is not complete, the priority vector cannot be calculated, so, in such case, the expert should be re-addressed with a request to provide the missing basic pair-wise comparisons, and, thus, ensure completeness of the set.
In the general case it is recommended to build spanning tree graphs using the unified vertex matrix (as prescribed in (Tsyganok, Kadenko & Andriichuk, 2015) and check *if at least one spanning tree* (basic pair-wise comparison set) can be built based on PCM, provided by all experts.
For instance, we can imagine the most extreme case when 3 experts are estimating 4 objects and the 1st expert provided only pair-wise comparison $a_{12}^{(1)}$, the 2nd expert – only pair-wise comparison $a_{23}^{(2)}$, and the 3rd expert – only pair-wise comparison $a_{34}^{(3)}$ (upper index denotes the number of the expert, who builds the PCM, while lower indices denote position of the pair-wise comparison in the PCM). Even in such extreme case the set of expert estimates can be considered complete, because we can build a spanning tree, covering all 4 objects for example, $O_1 \succ O_2 \succ O_3 \succ O_4$.
Focusing on extreme cases is not the purpose of this paper, as it will not allow us to illustrate the combinatorial approach that we are describing (these cases are not representative). So, further on we are going to consider cases when expert information is redundant and individual PCM, provided by every expert, allow us to build spanning trees:
$\forall k \in \{1..m\} T_k \neq 0$ where $k$ is the number of a particular expert and $T_k$ is the number of spanning trees, which can be built based on this expert's individual PCM. According to Caley's theorem on trees (Caley, 1889), $T_k \in (0,..,n^{n-2}]$.
In the context of individual PCM completeness we should mention a valuable and extensive research done by (Bozoki, Fulop, and Ronyai, 2010), providing optimal PCM completion rules.

### 4.3.2 Agreement

One of the innovative features of combinatorial algorithm is that it allows us to check and improve both consistency (i.e., “inner” consistency of the elements of a PCM, provided by an individual expert) and compatibility (i.e., “mutual” agreement of the respective elements of PCMs, provided by different experts) of pair-wise comparisons. In contrast to methods of aggregation of individual judgments (AIJ) and aggregation of individual priorities (AIP) (described, for instance, by (Forman & Peniwati, 1998), it is a one-phase method, allowing us to organize a step-by-step consistency/compatibility improvement procedure (as explained in (Tsyganok, Kadenko & Andriichuk, 2015).

The issue of relevance and necessity of consistency and compatibility improvement has been a subject of debate. However, in the class of weakly-structured problems, which we are dealing with, we assume that there exists some “ground truth” (the real ratios of interrelation between alternatives), so consistency and compatibility improvement *is necessary*.

In (Tsyganok, Kadenko & Andriichuk, 2015) it was recommended to use spectral consistency coefficient $K_c$, suggested in (Zgurovsky, Totsenko, & Tsyganok, 2004) as agreement measure. However, in (Olenko & Tsyganok, 2016) it was shown that this coefficient had its flaws. For example, originally $K_c$ was intended to lie within the range between 0 and 1, while in fact it could assume negative values (as shown by (Olenko & Tsyganok, 2016)). Beside that, introduction of $K_c$ was an attempt to unite two *independent* agreement measures: entropy and dispersion. This resulted in non-monotonous behavior of the coefficient. I.e., sometimes after pair-wise comparisons, lying on the largest distance from the average, were moved towards this average $K_c$ increased. In order to eliminate the drawbacks of the spectral consistency coefficient, (Olenko & Tsyganok, 2016) introduced the so-called double-entropy agreement index.

Let $N \geq 2$ denote the number of grades in the scale; $m \geq 2$ – the number of experts in the group; $r_i$ – the number of estimates that equal $i$ (i.e., the number of experts, who provided this estimate); $I = \{i_1, \ldots, i_k\}$– the set of expert estimates; $R = (r_1, \ldots, r_N)$– spectrum of expert estimates. In this agreement index two kinds of entropy (in the sense of “disagreement”) are used (this is Shannon entropy calculated by simplified Hartley formula (Hartley, 1928)):

**1-st**: $H(P) = -\sum_{j=1}^{k} p_j \ln p_j$ – among $k$ estimates;

**2-nd**: $H(Q) = -\sum_{i \in I} q_i \ln q_i$ – frequencies of $r_i$.

$P = \{p_1, \ldots, p_k\}$ are normalized distances between consecutive spectrum components: $p_j = d_j / d$, $j = 1, \ldots, k$. I.e. if $m$ experts provide their estimates in a scale with $n$ grades, their estimates (“votes”) are distributed among $k$ scale grades $\{i_1, \ldots, i_k\}$ $(2 \leq k \leq N)$. Distances are calculated as $d_j = i_{j+1} - i_j$, $j = 1, \ldots, k-1$, while the sum of these distances (by which they are divided for normalization) does not depend on specific scale values and always equals $d = \begin{cases} N - 1 + \left[\frac{N-1}{k-1}\right], & \text{if } k > 1 \\ N, & \text{if } k = 1 \end{cases}$ (as shown by (Olenko & Tsyganok, 2016)). Consequently, after normalization we get the “probabilities of appearance” or “shares” $P = \{p_1, \ldots, p_k\}$ of all distance values, which can be used in Shannon’s entropy formula in the 1st expression above.

In the 2nd expression $Q = \{q_i, i \in I\}$ are the normalized numbers of raters (experts) $\{r_1, \ldots, r_k\}$, who chose the respective scale values $I = \{i_1, \ldots, i_k\}$ as their estimates: $q_i = \frac{r_i}{\sum_{i \in I} r_i}$, $i \in I$.

After normalizing we obtain: $H^*(P)=\begin{cases}\dfrac{H(P)-\min_P H(P)}{\max_P H(P)-\min_P H(P)}, & \text{if } k<N\\ 1, & \text{if } k=N\end{cases}$ and $H^*(Q)=\dfrac{H(Q)}{\max_{k,Q} H(Q)}$. $\kappa_0(P,Q)=\dfrac{H^*(P)+H^*(Q)}{2}$ – disagreement index. And $\kappa(P,Q)=1-\kappa_0(P,Q)$ – agreement index, where $\kappa(P,Q)\in[0;1]$.

Presently, we suggest using this double entropy index instead of spectral consistency coefficient as agreement measure.

In order to calculate the double entropy index, spectrums should be built across pair-wise comparisons, provided by all the experts. In order to build spectrums, we suggest using the aforementioned combinatorial approach. Let us describe consistency verification and improvement procedures in greater detail. We assume that at this stage of the process, completeness of expert data is already ensured. As completeness is reached, we can use the following step-by-step consistency improvement procedure:

Step 1: enumerate all basic sets of pair-wise comparisons (spanning trees) in each individual pair-wise comparison matrix.

If dimensionality of the PCM is *n*x*n* and this PCM is complete (there are no missing elements), the number of spanning trees will amount to $n^{n-2}$, according to the aforementioned Caley's "theorem on trees" (Caley, 1889). Otherwise (if PCM are incomplete and some comparisons are missing), we can assume that individual PCM, provided by expert number *k* allows us to build $T_k$ spanning trees. The total number of trees *T* is limited from above: $T=\sum_{k=1}^{m} T_k \le mn^{n-2}$.

Step 2: build ideally consistent pair-wise comparison matrices.

Each of basic sets of pair-wise comparisons allows us to "reconstruct" an ideally consistent pair-wise comparison matrix (ICPCM), following the transitivity rule. For example, if a spanning tree covers alternatives $A_l$, $A_i$, and $A_j$, then $a_{ij}^{reconstructed}=a_{il}\times a_{lj}$. As a result we will get *T* ICPCMs.

Step 3: calculate priorities based on each ICPCM.

As matrices are ideally consistent, priorities can be calculated as normalized elements of any basic pair-wise comparison set, for example, the elements from any ICPCM line.

$$\begin{pmatrix} a_{11} & & a_{1n} \\ & \cdots & \\ a_{n1} & & a_{nn}\end{pmatrix}=\begin{pmatrix} 1 & & \frac{w_1}{w_n} \\ & \cdots & \\ \frac{w_n}{w_1} & & 1\end{pmatrix}\Rightarrow \forall i,j=1..n: w_j=\frac{a_{ij}}{\sum_{l=1}^{n} a_{il}}$$

As a result of this step we will get *T* priority vectors: $\{w^q=\{w_j^q : j=1..n\};q=1..T\}$.

Step 4: Build priority spectrums.

Spectrums are built for each coordinate of priority vectors (as shown in (Zgurovsky, Totsenko, & Tsyganok, 2004) and (Olenko & Tsyganok, 2016)). Values are rounded to pre-defined precision threshold (say, 0.01 or 0.001). As a result we get *n* spectrums, each consisting of *T* priority estimates. If experts have different competencies, the estimates should be multiplied by respective experts' relative competence coefficients (see problem statement) before being included into the spectrums. Typical looks of spectrums can be seen on Figures 2a – 2d in Section 5 of the current paper.

Step 5: Calculate double entropy coefficients and compare them with thresholds.

Double entropy coefficients (indices) are calculated as explained in the beginning of this sub-section. As a result we get *n* double entropy index values (one for each priority vector coordinate). (Olenko & Tsyganok 2016) choose the value of 0.7 as usability threshold (instead of usability threshold, suggested in (Zgurovsky, Totsenko, & Tsyganok, 2004)). Empirical considerations behind the choice of this particular value are provided in their publication however we have to admit that development of an analytical expression double entropy threshold should be a subject of a separate research, lying beyond the scope of this paper.

If all $n$ double entropy indices ( $K_l; l = 1..n$ ) are above the threshold value ( $\min_{l=1}^{n} K_l > 0.7$ ), then priority vectors can be aggregated, and aggregate priority vector $w^{aggregate}$ can be accepted as the final result (solution of the problem, posed in section 3).
Otherwise, agreement improvement is necessary. In order to improve agreement of expert estimates, we have to "move" the most "outstanding" of respective pair-wise comparisons towards the average, i.e. request the respective experts to change their respective initial judgments. In order to define this average, we suggest building an aggregate ICPCM, $A^{aggregate} = \{a_{ij}^{aggregate}; i.j = 1..n\}$ whose elements will represent the "target", towards which we are going to move in the process of consistency improvement. Such aggregate ICPCM can be easily reconstructed based on the aggregate priority vector.

$$a_{ij}^{aggregate} = \frac{w_i^{aggregate}}{w_j^{aggregate}}; i, j = 1..n \quad (3)$$

As we can see, we've arrived at a point when the next step (in any of the two listed cases) is the aggregation of priority values. This brings us to the next phase of problem solution process.

### 4.4 Aggregation of estimates

In principle, the aggregation procedure is based on calculation of geometric (or arithmetic) mean across all $T$ priority vectors, obtained at the previous stage of problem solution. The simplest formula for calculation of aggregate vector would be:

$$w_j^{aggregate} = (\prod_{q=1}^{T} w_j^q)^{\frac{1}{T}}; j = 1..n \quad (4)$$

However, this formula does not take any aspects of expert competence and confidence into account. (Tsyganok, Kadenko, & Andriichuk, 2012) showed that individual expert competence can be neglected if the number of experts in the group is above several decades, while in small groups, the differences in individual experts' relative competence should be taken into consideration.

We feel, that in the context of our problem, expert competence can (and should) be considered at several conceptual levels:

1) Subject domain and specific problem
2) Specific matrix
3) Specific spanning tree
4) Specific pair-wise comparison

Estimation of relative expert competence levels in a group in the context of subject domain in general and specific problem is a very complex task, which is not the subject of current research. We should only note that expert competence at this highest conceptual level can be based on self-estimation, mutual estimation, and objective component (reflecting the overall experience and achievements of the expert). Beside that, specific problem (or strategic goal) can be described by basic keywords, and relative expert competence levels can be defined in relation to each specific keyword, describing the problem (or overall goal). In this paper we can assume that relative competence levels, given in the problem statement in section 3 of this paper reflect the first (most general) conceptual level.

When it comes to a particular individual PCM, we feel that the expert competence level should be reflected by PCM completeness, its consistency, and compatibility with individual PCMs of other experts.

In terms of completeness, we should say, that the "more complete" the individual PCM built by expert number $k$ is, the more spanning trees $T_k$ can be reconstructed on its basis. Consequently, the relative number of multipliers in formula (4) ($T_k/T$) (reflecting the significance of the respective expert's estimates) will increase.

In order to consider consistency and compatibility while aggregation of priorities, we suggest assigning a rating to each ICPCM (and the respective priority vector), built based on an individual PCM of a specific expert. This rating should reflect the inner agreement of pair-wise comparisons from the PCM of this expert (i.e., consistency) and the mutual agreement between pair-wise comparisons PCMs of different experts (i.e., compatibility). For this purpose, we suggest building $m$ replicas (copies) of each ICPCM. As a result we will get $T^{*}=mT$ ICPCM ($m \le T^{*} \le m^2 n^{n-2}$). Beside consistency and compatibility, the rating should also incorporate expert competence level in terms of scale (s)he uses (the more grades in the scale used, the more competence/confidence) and relative competence levels, set a-priori in the problem statement in section 3. Based on these considerations, we suggest building a rating of each specific ICPCM (and the respective priority vector) out of these $T^*$ matrices according to the following formula:

$$R_{kql} = \frac{c_k c_l s^{kq} s^l}{\ln\left(\sum_{u,v} \left|a_{uv}^{kq} - a_{uv}^{l}\right| + e\right)} \quad (5)$$

where $k$, $l$ are the numbers of two experts ($k,l=1..m$), whose matrices are compared; $c_k, c_l$ are relative competence coefficients from the problem statement in section 3 of this paper; $k$ and $l$ can be equal or different, i.e. ICPCM, based on the individual PCM of a given expert number $k$ can be compared with this individual PCM and with individual PCMs built by other experts; $q$ is the number of the respective ICPCM copy ($q=1..mT_k$); $s^{kq}$ is the average weight of scales, in which the pair-wise comparisons from the basic set (spanning tree) number $q$, reconstructed based on PCM of expert number $k$, were input ($s^{kq}$ is calculated according to formula (1) above); $s^l$ is the average weight of scales, in which the elements of PCM of expert number $l$ were built ($s^l$ is calculated based on formula (2) above).
If we incorporate normalized ICPCM rating values into formula (4), we get the following expression for aggregate priorities:

$$w_j^{aggregate} = \prod_{k,l=1}^{m} \left( \prod_{q^k=1}^{T_k} \left( w_j^{(kq_k l)} \right)^{\frac{R_{kq_k l}}{\sum_{u,p,v} R_{upv}}} \right); j = 1..n \quad (6)$$

When aggregate priorities are calculated, we can either accept their values as the final result and the problem solution, or improve the level of estimates' agreement.

### 4.5 Agreement level improvement (if necessary)

Agreement improvement is necessary if double entropy index for some of priority vector coordinates lies below the threshold value: $\exists l = 1..n : K_l \le 0.7$. In order to improve agreement of expert estimates (as mentioned above), we should move the most inconsistent and incompatible pair-wise comparisons towards respective average values. For this purpose, we should request the respective experts to change their respective pair-wise comparisons accordingly.
The average values are the elements of aggregate ICPCM (see formula (3)). The element, that should be changed in order to improve the double entropy index for priority number $l$ is

$$a_{lj*}^{(k)} = \arg\max_{k=1..m;\, j=1..n} \left| a_{lj}^{aggregate} - a_{lj}^{(k)} \right| \quad (7)$$

After each feedback step aggregate ICPCM is recalculated. The process continues until the agreement condition is satisfied: $\min_{l=1}^{n} K_l > 0.7$.

## 5. Numeric example

In this section we consider a numeric example, where 3 equally competent experts (from left to right) provide pair-wise comparisons of a set of 4 alternatives. If an expert "skips" some pair-wise comparison, the respective cell in the PCM remains blank and respective spanning trees are not built. Numeric data for the example is shown on Figures 1a and 1b. Spectrums of aggregate priorities $w_1,...,w_4$, rounded to the nearest scale grades {1,..,9} are shown on Figures 2a – 2d. On these figures scale grades are shown along axis $X$, while the ratings (numbers of IC PCMs, from which the respective weight values are derived) are shown along axis $Y$.
For the sake of simplicity in the example we consider that all experts provide estimates in the integer scale and as all double entropy values are above the threshold. Normalized average priority vector values are: $\{w_1,...,w_4\}$ = {0.0918; 0.1908; 0.1808; 0.5366}.

**ORDINAL PREFERENCES**

| 1 | 2 | 3 | 4 | 1 | 2 | 3 | 4 | 1 | 2 | 3 | 4 |
|---|---|---|---|---|---|---|---|---|---|---|---|
| 1 | < | * | < | 1 | < | > | * | 1 | < | > | < |
|  | 1 | > | < |  | 1 | > | < |  | 1 | < | * |
|  |  | 1 | < |  |  | 1 | < |  |  | 1 | < |
|  |  |  | 1 |  |  |  | 1 |  |  |  | 1 |

**NUMBER OF GRADES**

| | | | | | | | | | | | |
|---|---|---|---|---|---|---|---|---|---|---|---|
| 1 | 7 | * | 9 | 1 | 9 | 7 | * | 1 | 9 | 9 | 9 |
|  | 1 | 5 | 9 |  | 1 | 7 | 7 |  | 1 | 9 | * |
|  |  | 1 | 9 |  |  | 1 | 5 |  |  | 1 | 9 |
|  |  |  | 1 |  |  |  | 1 |  |  |  | 1 |

**SCALE GRADE NUMBER**

| | | | | | | | | | | | |
|---|---|---|---|---|---|---|---|---|---|---|---|
| 1 | 2 | * | 8 | 1 | 2 | 4 | * | 1 | 2 | 4 | 8 |
|  | 1 | 2 | 4 |  | 1 | 2 | 4 |  | 1 | 2 | * |
|  |  | 1 | 2 |  |  | 1 | 2 |  |  | 1 | 2 |
|  |  |  | 1 |  |  |  | 1 |  |  |  | 1 |

**INDIVIDUAL PAIR-WISE COMPARISON MATRICES OF**

**EXPERT 1**

| 1 | 2 | 3 | 4 |
|---|---|---|---|
| 1 | 1/2 | 1 | 1/8 |
| 2 1/6 | 1 | 2 1/2 | 1/4 |
| 1 | 2/5 | 1 | 1/2 |
| 8 | 4 | 2 | 1 |

**EXPERT 2**

| 1 | 2 | 3 | 4 |
|---|---|---|---|
| 1 | 1/2 | 4 5/6 | 1 |
| 2 | 1 | 2 1/6 | 1/5 |
| 1/5 | 1/2 | 1 | 2/5 |
| 1 | 4 5/6 | 2 1/2 | 1 |

**EXPERT 3**

| 1 | 2 | 3 | 4 |
|---|---|---|---|
| 1 | 1/2 | 4 | 1/8 |
| 2 | 1 | 1/2 | 1 |
| 1/4 | 2 | 1 | 1/2 |
| 8 | 1 | 2 | 1 |

**SPANNING TREES**

**IDEALLY CONSISTENT PAIR-WISE COMPARISON MATRICES**

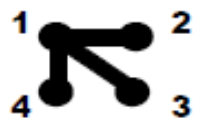

| | | | | | | | | | | | |
|---|---|---|---|---|---|---|---|---|---|---|---|
|  |  |  |  |  |  |  |  | 1 | 1/2 | 4 | 1/8 |
|  |  |  |  |  |  |  |  | 2 | 1 | 8 | 1/4 |
|  |  |  |  |  |  |  |  | 1/4 | 1/8 | 1 | 0 |
|  |  |  |  |  |  |  |  | 8 | 4 | 32 | 1 |

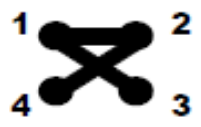

| | | | | | | | | | | | |
|---|---|---|---|---|---|---|---|---|---|---|---|
|  |  |  |  | 1 | 1/2 | 4 5/6 | 1/9 |  |  |  |  |
|  |  |  |  | 2 | 1 | 9 2/3 | 1/5 |  |  |  |  |
|  |  |  |  | 1/5 | 1/9 | 1 | 0 |  |  |  |  |
|  |  |  |  | 9 2/3 | 4 5/6 | 46 5/7 | 1 |  |  |  |  |

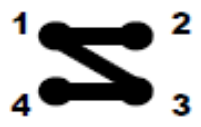

| | | | | | | | | | | | |
|---|---|---|---|---|---|---|---|---|---|---|---|
|  |  |  |  | 1 | 1/2 | 4 5/6 | 2 | 1 | 1/2 | 4 | 2 |
|  |  |  |  | 2 | 1 | 9 2/3 | 3 7/8 | 2 | 1 | 8 | 4 |
|  |  |  |  | 1/5 | 1/9 | 1 | 2/5 | 1/4 | 1/8 | 1 | 1/2 |
|  |  |  |  | 1/2 | 1/4 | 2 1/2 | 1 | 1/2 | 1/4 | 2 | 1 |

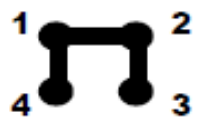

| | | | | | | | | | | | |
|---|---|---|---|---|---|---|---|---|---|---|---|
| 1 | 1/2 | 1 1/7 | 1/8 |  |  |  |  | 1 | 1/2 | 1/4 | 1/8 |
| 2 1/6 | 1 | 2 1/2 | 1/4 |  |  |  |  | 2 | 1 | 1/2 | 1/4 |
| 7/8 | 2/5 | 1 | 1/9 |  |  |  |  | 4 | 2 | 1 | 1/2 |
| 8 | 3 2/3 | 9 2/9 | 1 |  |  |  |  | 8 | 4 | 2 | 1 |

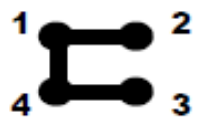

| | | | | | | | | | | | |
|---|---|---|---|---|---|---|---|---|---|---|---|
| 1 | 1/2 | 1/4 | 1/8 |  |  |  |  | 1 | 1/2 | 1/4 | 1/8 |
| 2 1/6 | 1 | 1/2 | 1/4 |  |  |  |  | 2 | 1 | 1/2 | 1/4 |
| 4 | 1 6/7 | 1 | 1/2 |  |  |  |  | 4 | 2 | 1 | 1/2 |
| 8 | 3 2/3 | 2 | 1 |  |  |  |  | 8 | 4 | 2 | 1 |

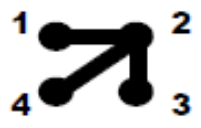

| | | | | | | | | | | | |
|---|---|---|---|---|---|---|---|---|---|---|---|
| 1 | 1/2 | 1 1/7 | 1/9 | 1 | 1/2 | 1 | 1/9 |  |  |  |  |
| 2 1/6 | 1 | 2 1/2 | 1/4 | 2 | 1 | 2 1/6 | 1/5 |  |  |  |  |
| 7/8 | 2/5 | 1 | 0 | 1 | 1/2 | 1 | 0 |  |  |  |  |
| 8 2/3 | 4 | 10 | 1 | 9 2/3 | 4 5/6 | 10 1/2 | 1 |  |  |  |  |

Fig. 1a Numeric example data

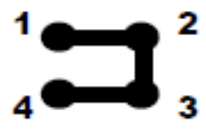

| | | | |
|---|---|---|---|
| 1 | 1/2 | 1 1/7 | 4/7 |
| 2 1/6 | 1 | 2 1/2 | 1 1/4 |
| 7/8 | 2/5 | 1 | 1/2 |
| 1 3/4 | 4/5 | 2 | 1 |
| 1 | 1/2 | 1 | 3/7 |
| 2 | 1 | 2 1/6 | 6/7 |
| 1 | 1/2 | 1 | 2/5 |
| 2 1/3 | 1 1/6 | 2 1/2 | 1 |
| 1 | 1/2 | 1/4 | 1/8 |
| 2 | 1 | 1/2 | 1/4 |
| 4 | 2 | 1 | 1/2 |
| 8 | 4 | 2 | 1 |

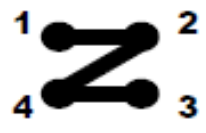

| | | | |
|---|---|---|---|
| 1 | 1/2 | 1/4 | 1/9 |
| 2 1/6 | 1 | 1/2 | 1/4 |
| 4 1/3 | 2 | 1 | 1/2 |
| 8 2/3 | 4 | 2 | 1 |
| 1 | 1/2 | 1/4 | 1/9 |
| 2 | 1 | 1/2 | 1/5 |
| 3 7/8 | 2 | 1 | 2/5 |
| 9 2/3 | 4 5/6 | 2 1/2 | 1 |

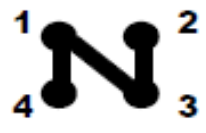

| | | | |
|---|---|---|---|
| 1 | 8 | 4 | 1/8 |
| 1/8 | 1 | 1/2 | 0 |
| 1/4 | 2 | 1 | 0 |
| 8 | 64 | 32 | 1 |

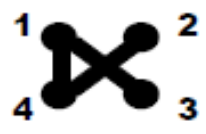

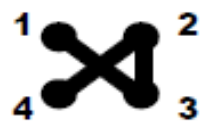

| | | | |
|---|---|---|---|
| 1 | 2 1/4 | 4 5/6 | 1/2 |
| 4/9 | 1 | 2 1/6 | 1/5 |
| 1/5 | 1/2 | 1 | 0 |
| 2 1/6 | 4 5/6 | 10 1/2 | 1 |

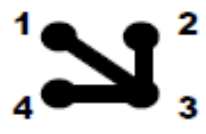

| | | | |
|---|---|---|---|
| 1 | 2 1/4 | 4 5/6 | 2 |
| 4/9 | 1 | 2 1/6 | 6/7 |
| 1/5 | 1/2 | 1 | 2/5 |
| 1/2 | 1 1/6 | 2 1/2 | 1 |
| 1 | 8 | 4 | 2 |
| 1/8 | 1 | 1/2 | 1/4 |
| 1/4 | 2 | 1 | 1/2 |
| 1/2 | 4 | 2 | 1 |

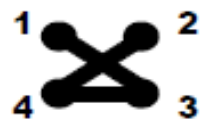

| | | | |
|---|---|---|---|
| 1 | 9 1/3 | 4 5/6 | 2 |
| 1/9 | 1 | 1/2 | 1/5 |
| 1/5 | 2 | 1 | 2/5 |
| 1/2 | 4 5/6 | 2 1/2 | 1 |

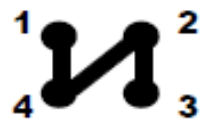

| | | | |
|---|---|---|---|
| 1 | 1/2 | 1 1/4 | 1/8 |
| 2 | 1 | 2 1/2 | 1/4 |
| 4/5 | 2/5 | 1 | 0 |
| 8 | 4 | 10 | 1 |

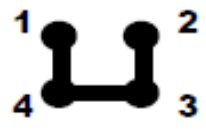

| | | | |
|---|---|---|---|
| 1 | 0 | 1/4 | 1/8 |
| 10 | 1 | 2 1/2 | 1 1/4 |
| 4 | 2/5 | 1 | 1/2 |
| 8 | 4/5 | 2 | 1 |
| 1 | 1/2 | 1/4 | 1/8 |
| 2 | 1 | 1/2 | 1/4 |
| 4 | 2 | 1 | 1/2 |
| 8 | 4 | 2 | 1 |

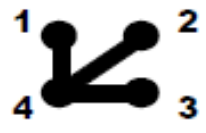

| | | | |
|---|---|---|---|
| 1 | 1/2 | 1/4 | 1/8 |
| 2 | 1 | 1/2 | 1/4 |
| 4 | 2 | 1 | 1/2 |
| 8 | 4 | 2 | 1 |

Fig. 1b Numeric example data (continued)

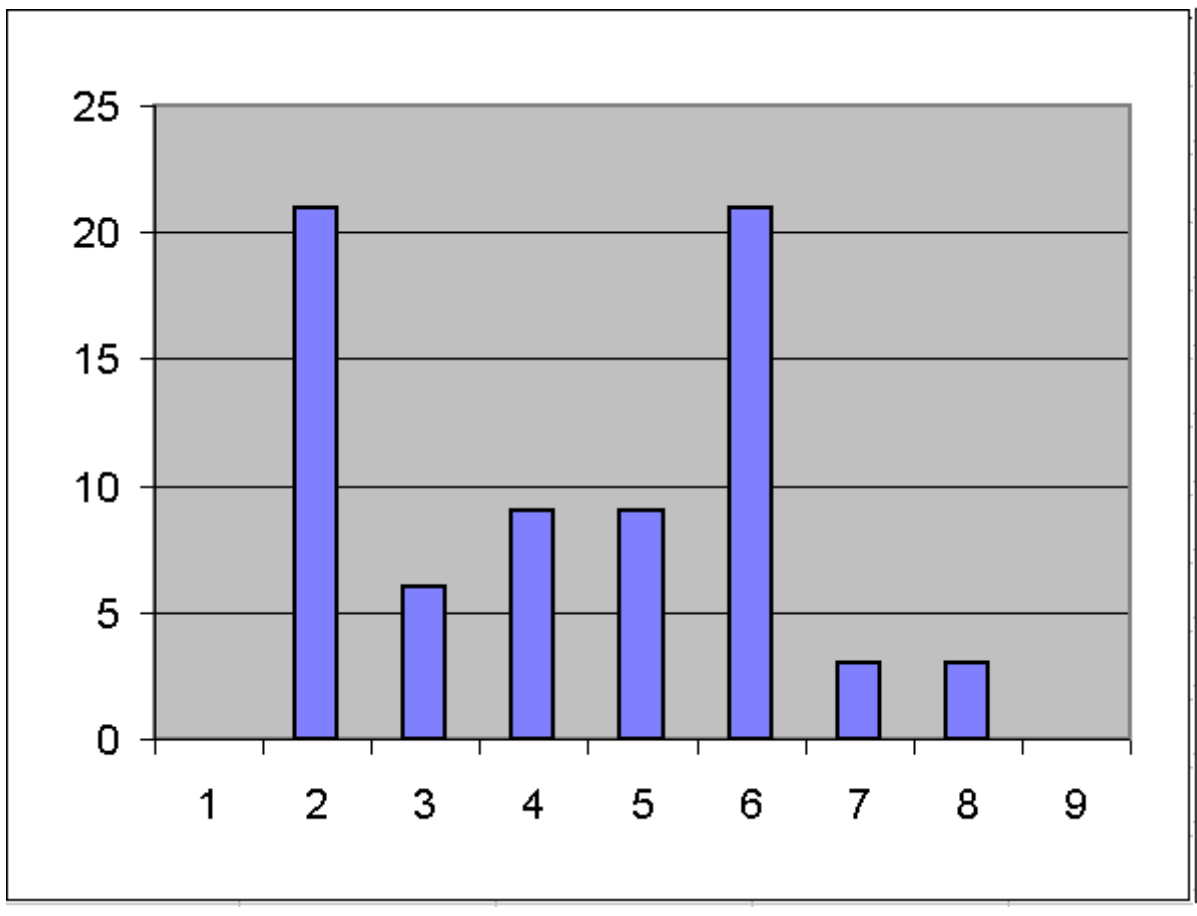

Figure 2a Spectrum of $w_1$

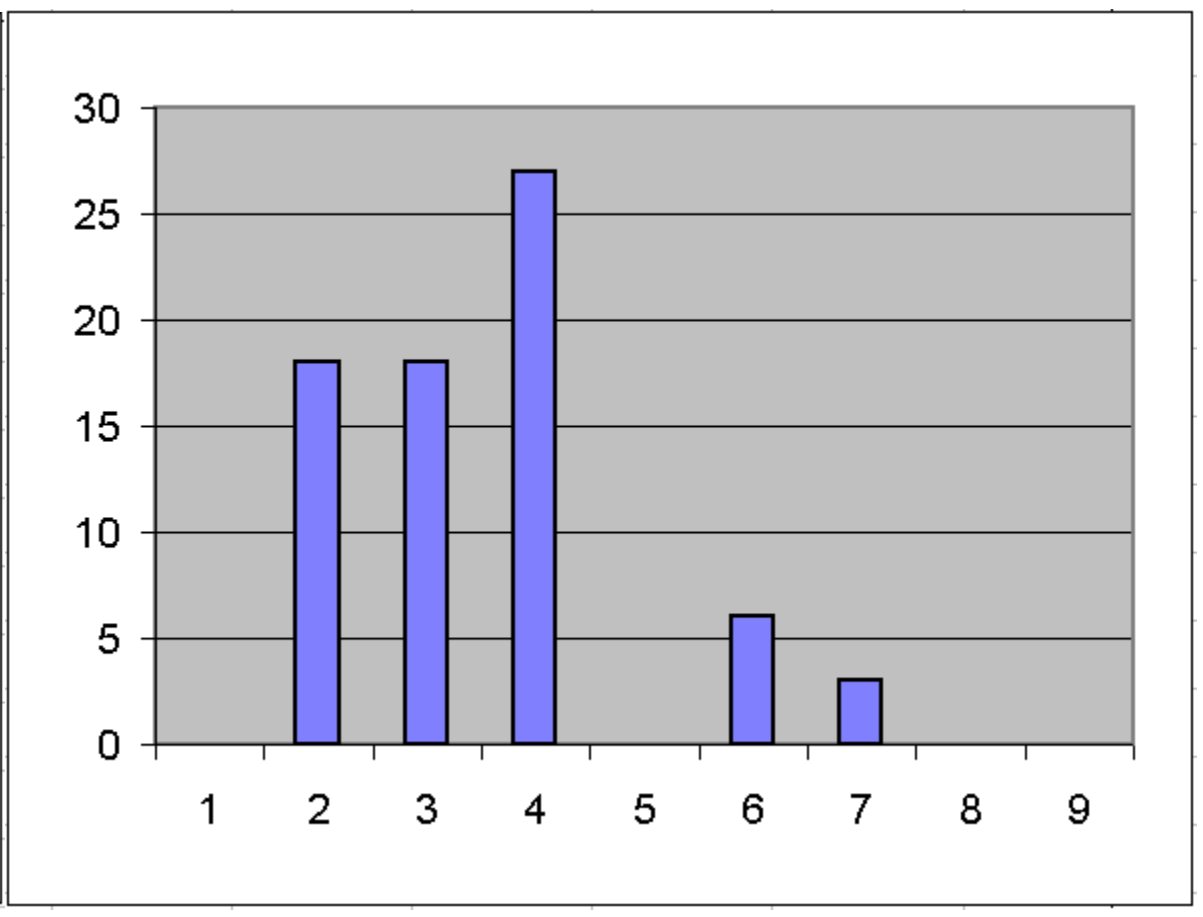

Figure 2b Spectrum of $w_2$

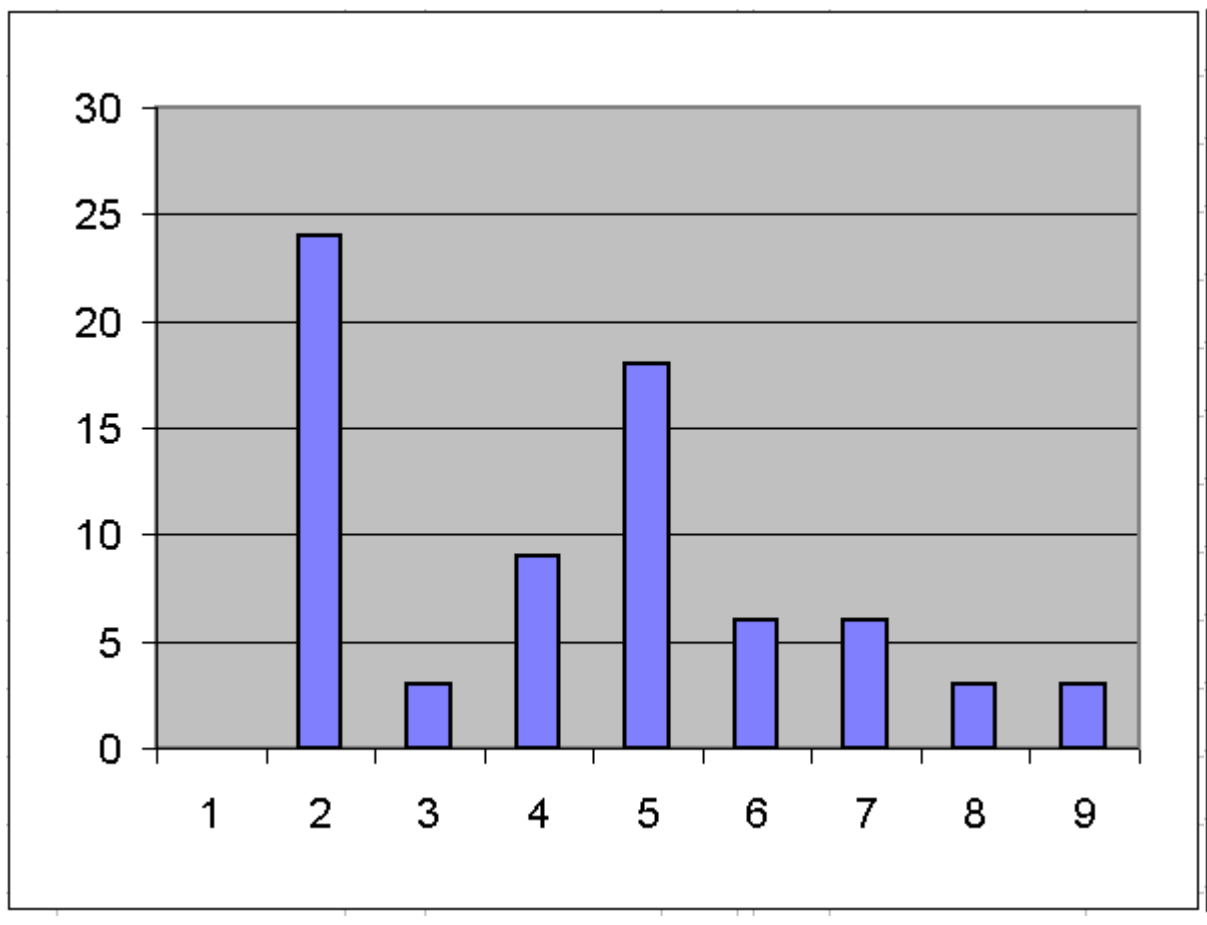

Figure 2c Spectrum of $w_3$

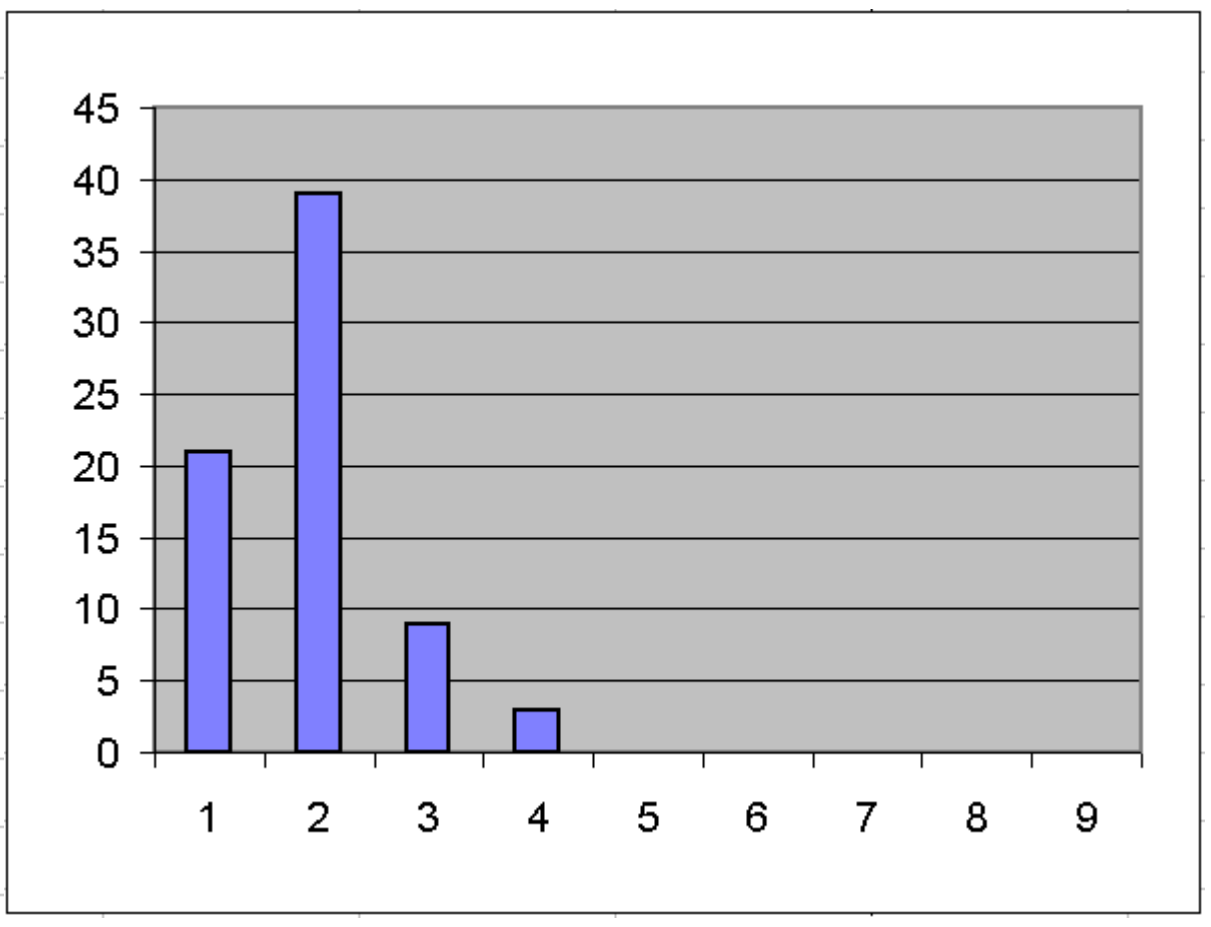

Figure 2d Spectrum of $w_4$

## 6. Conclusions

A Method with Feedback for Aggregation of Group Incomplete Pair-Wise Comparisons using Scales with Different Numbers of Grades is suggested. The method allows the expert group to reach consensus, which is necessary for aggregation of individual estimates. The advantages of the method are as follows. The method

1) utilizes redundancy of expert information most thoroughly
2) allows to aggregate pair-wise comparisons, obtained from a group of experts
3) is suitable for both complete and incomplete PCMs
4) allows to organize step-by-step feedback with individual experts
5) allows to improve both consistency and compatibility of individual PCMs
6) can use both geometric mean and arithmetic mean as the principal aggregation rule
7) allows experts to use different scales for each pair-wise comparison
8) allows to take all the aspects of expert competence into consideration.

While the mathematical aspect of the method may seem over-complicated, we should stress, that all the calculations are to be performed by DSS software. The experts only need to input pair-wise comparisons in the scale which they find most suitable.

## 7. Key References